\documentclass[preprint,11pt,times]{elsarticle}

\usepackage{amsmath,amssymb,amsfonts,amsthm}
\usepackage{graphicx,multirow,mathrsfs}
\usepackage[title]{appendix}
\usepackage{xcolor,textcomp,manyfoot,booktabs}
\usepackage{algorithm,algorithmicx,algpseudocode,listings}
\usepackage{geometry}
\usepackage{setspace}
\usepackage[
colorlinks=true,
linkcolor=blue,
anchorcolor=blue,
citecolor=blue
]{hyperref}

\newcommand{\manuscriptheadfont}{%
	\normalfont\fontseries{bx}\selectfont\boldmath}
\newcommand{\manuscriptsubheadfont}{%
	\normalfont\fontseries{bx}\itshape\selectfont\boldmath}

\newtheoremstyle{manuscriptplain}
{\topsep}{\topsep}{\itshape}{}
{\manuscriptheadfont}{.}{0.5em}
{\thmname{#1}\thmnumber{ #2}\thmnote{ (#3)}}

\newtheoremstyle{manuscriptupright}
{\topsep}{\topsep}{\normalfont}{}
{\manuscriptheadfont}{.}{0.5em}
{\thmname{#1}\thmnumber{ #2}\thmnote{ (#3)}}

\theoremstyle{manuscriptplain}
\newtheorem{theorem}{Theorem}[section]
\newtheorem{corollary}[theorem]{Corollary}
\newtheorem{lemma}[theorem]{Lemma}
\newtheorem{proposition}[theorem]{Proposition}

\theoremstyle{manuscriptupright}

\newtheorem{remark}{Remark}[section]
\newtheorem{example}{Example}[section]
\numberwithin{equation}{section}

\makeatletter
\AtBeginDocument{%
	\renewcommand{\section}{%
		\@startsection{section}{1}{\z@}%
		{18\p@ \@plus 6\p@ \@minus 3\p@}%
		{9\p@ \@plus 6\p@ \@minus 3\p@}%
		{\normalsize\manuscriptheadfont}}%
	\renewcommand{\subsection}{%
		\@startsection{subsection}{2}{\z@}%
		{12\p@ \@plus 6\p@ \@minus 3\p@}%
		{3\p@ \@plus 6\p@ \@minus 3\p@}%
		{\normalsize\manuscriptsubheadfont}}%
	\renewcommand{\subsubsection}{%
		\@startsection{subsubsection}{3}{\z@}%
		{12\p@ \@plus 6\p@ \@minus 3\p@}%
		{\p@}%
		{\normalsize\manuscriptsubheadfont}}%
}
\makeatother

\providecommand{\backmatter}{}
\providecommand{\bmhead}[1]{%
	\par\medskip\noindent{\manuscriptheadfont #1}\quad}

\begin{document}
	
	\begin{frontmatter}
		
		\title{ Hypercyclic Algebras for Weighted Double Shifts}
		
		\author[aff1]{Xiang Chen\corref{cor1}}
		\ead{2020233018@tju.edu.cn}

		\author[aff2]{Wei Wang}
		\ead{wangwei223@htu.edu.cn}
		
		\author[aff1]{Ze-hua Zhou}
		\ead{zhzhou@tju.edu.cn}
		
		\cortext[cor1]{Corresponding author}
		
		\affiliation[aff1]{organization={School of Mathematics, Tianjin University},
			city={Tianjin}, postcode={300350},
			country={P.R. China}}
		
		\affiliation[aff2]{organization={School of Mathematics and Statistics, Henan Normal University},
			city={Xinxiang}, postcode={453007},
			country={P.R. China}}
	\begin{abstract}
		Building on the approach of Bayart, Costa J\'unior and
		Papathanasiou, we establish a Baire criterion for possibly
		noncommuting operator pairs to support hypercyclic algebras,
		with hypercyclicity defined through all finite compositions.
		We then study weighted double backward shifts under
		different multiplications and topologies.
		We characterize the existence of hypercyclic algebras for
		coordinatewise multiplication on sequence algebras admitting
		a continuous norm and, separately, on
		$\omega^2=\mathbb C^{\mathbb N_0^2}$ with the product topology.
		For the noncommutative triangular product induced by upper
		triangular matrix multiplication, we obtain necessary and
		sufficient conditions for the existence of dense hypercyclic
		algebras on $\omega^2$.
		In this setting, their existence is equivalent to that of
		one freely generated by countably infinitely many elements.
		Finally, we exhibit a pair on a fixed Fr\'echet space that
		supports a dense hypercyclic algebra for coordinatewise
		multiplication and a nonzero hypercyclic algebra, but no
		dense one, for the triangular product.
	\end{abstract}
		
		
	\end{frontmatter}
	
	\section{Introduction}
	
	Let $X$ be a separable Fr\'echet space,  and let
	$T:X\to X$ be a continuous linear operator.
	The set of hypercyclic vectors for $T$ is
	\[
	HC(T)
	=
	\left\{
	x\in X:
	\overline{\{T^n x:n\in\mathbb N_0\}}=X
	\right\}.
	\]
	The operator $T$ is called \emph{hypercyclic} if
	$HC(T)\ne\varnothing$.
	Weighted backward shifts provide basic models for studying
	this property; see~\cite{BM,CC} for general background.
	Suppose, in addition, that $X$ is an associative Fr\'echet
	algebra with jointly continuous multiplication $\ast$.
	A nonzero subalgebra $\mathcal A\subseteq X$ is called a
	\emph{hypercyclic algebra} for $T$ if
	\[
	\mathcal A\setminus\{0\}\subseteq HC(T).
	\]
	The subalgebra need not be closed or contain an identity,
	and the multiplication need not be commutative.
	
	We equip the complex double sequence space
	\[
	\omega^2
	=
	\mathbb C^{\mathbb N_0\times\mathbb N_0}
	=
	\{x=(x_{i,j})_{i,j\ge0}:x_{i,j}\in\mathbb C\}
	\]
	with the product topology.
	The topology and convergence properties of double sequence
	spaces have been studied
	in~\cite{AB2005,BS2009,BS1999,ME2003,TMRY2025}.
	Approximation in the product topology involves finitely many
	coordinates at a time, while summability norms also require
	control of the tails.
	Shu, Wang and Zhao~\cite{SWZ2017} characterized hypercyclicity,
	weak mixing and mixing for weighted double backward shifts.
	Wang~\cite{Wang2026} studied hypercyclicity, mixing and chaos
	for unweighted double shifts on double Fr\'echet sequence spaces.
	With variable weights, the coordinate shifts need not commute,
	and compositions with the same displacement may have
	different weight products.
	
B\`es and Papathanasiou~\cite{BP2020} constructed dense
hypercyclic algebras that are not finitely generated for
convolution operators.
For the Cauchy product on $\ell^1$, Falc\'o and
Grosse-Erdmann~\cite{FG2020Shifts} proved that mixing weighted
backward shifts have algebrable sets of hypercyclic vectors;
in particular, their result applies to the Rolewicz operators.
They also obtained algebrability for the MacLane
differentiation operator on $H(\mathbb C)$ with its usual
pointwise multiplication.
Bayart, Costa J\'unior and Papathanasiou~\cite{BCP2021}
developed a Baire method based on separating finitely many
monomials, applicable also to noncommutative algebras,
and characterized weighted shifts supporting hypercyclic
algebras on Fr\'echet sequence algebras.
Further work concerns disjoint and common hypercyclic
algebras~\cite{BCP2022}, frequently hypercyclic
algebras~\cite{FGF2020}, weighted shifts on
trees~\cite{AC2026,CZZ2026}, and discrete composition
operators~\cite{ZCZ2026}.
Recently, Costa Jr. and Rocha~\cite{CR2026} proved that every
separable infinite-dimensional Fr\'echet algebra admits a
continuous linear operator supporting a dense invariant
hypercyclic algebra.
This existence theorem does not provide criteria for
a prescribed pair of weighted shifts.
	
	Here we study pairs $\mathcal T=(T_1,T_2)$ of weighted
	double backward shifts, without assuming commutativity.
	Hypercyclicity refers to the orbit under all finite
	compositions, including the identity.
	This notion differs from common or disjoint hypercyclicity.
	Finite compositions are described by paths, whose weight
	products may depend on the order of the factors.
	
	Our main concern is the dependence of dense hypercyclic
	algebras on the multiplication.
	Besides coordinatewise multiplication, we consider the
	associative triangular product
	\[
	(x\ast y)_{i,j}
	=
	\sum_{t=0}^{j}x_{i,t}y_{i+t,j-t}.
	\]
	Under the identification $x_{i,j}\leftrightarrow A_{i,i+j}$,
	this is upper triangular matrix multiplication.
	It couples distinct coordinates and preserves word order.
	With this product, $\omega^2$ is a Fr\'echet algebra in the
	product topology, and $\mathcal L_p=\ell^p(\mathbb N_0^2)$
	is a Banach algebra in its usual norm for $1\le p\le2$.
	
Taking the approach of~\cite{BCP2021} as a starting point,
we study hypercyclic algebras for possibly noncommuting
pairs of weighted double backward shifts on double
Fr\'echet sequence algebras $X\subseteq\omega^2$.
Our focus is how the multiplication and the topology
affect the existence, density and free generation of
such algebras.
We consider coordinatewise and triangular multiplication,
together with the product topology on $\omega^2$ and
the usual norm topologies on summable sequence spaces.
For coordinatewise sequence algebras admitting a continuous
norm, we obtain a characterization using inverse powers
of path weights and translated coordinate vectors;
$\omega^2$ is treated separately through compact families
of polynomial coefficients.
For triangular multiplication, we establish path-weight
characterizations on $\omega^2$ and $\mathcal L_p$,
$1\le p\le2$, distinguishing nonvanishing weights in
the product topology from uniform growth in norm.
In both settings, the existence of a dense hypercyclic
algebra is equivalent to that of one freely generated
by countably infinitely many elements.
The restriction $p\le2$ is sharp for continuity of this
product in the usual norm.
Finally, an example shows that, with the space, topology
and operator pair fixed, changing only the multiplication
can preserve the existence of nonzero hypercyclic algebras
while eliminating dense ones.

Section~\ref{sec:baire} establishes the general Baire criterion.
Section~\ref{sec:pointwise} treats coordinatewise multiplication.
Section~\ref{sec:triangular-product} develops the
triangular-product characterizations and concludes with
the example comparing the two multiplications.

	\section{A Baire criterion for operator pairs}\label{sec:baire}
	
	Throughout this section, $X$ is a nonzero separable associative
	Fr\'echet algebra over $\mathbb C$,
	with jointly continuous multiplication $\ast$.
	Let $\mathcal T=(T_1,T_2)$ be a pair of continuous linear
	operators on $X$.
	Neither the algebra nor the operators are assumed to commute.
	
	Let
	\[
	\mathcal F_{\mathcal T}
	=
	\{I\}\cup
	\{T_{i_1}\cdots T_{i_r}:
	r\ge1,\ i_1,\ldots,i_r\in\{1,2\}\}.
	\]
	A \emph{finite path} is a word $\sigma=(i_1,\ldots,i_r)$
	over $\{1,2\}$, with associated operator
	$T^\sigma=T_{i_1}\cdots T_{i_r}$.
	The empty path corresponds to $I$.
	Define
	\[
	\begin{aligned}
		\operatorname{orb}(\mathcal T,x)
		=\{Sx:S\in\mathcal F_{\mathcal T}\},\quad
		HC(\mathcal T)=\left\{
		x\in X:
		\overline{\operatorname{orb}(\mathcal T,x)}=X
		\right\}.
	\end{aligned}
	\]
	The pair $\mathcal T$ is hypercyclic if
	$HC(\mathcal T)\ne\varnothing$.
	A nonzero subalgebra $\mathcal A\subseteq X$ is a hypercyclic
	algebra for $\mathcal T$ if
	$\mathcal A\setminus\{0\}\subseteq HC(\mathcal T)$.
	
	For $d\ge1$, let $\mathscr P_d$ consist of the formal polynomials
	\[
	P(z_1,\ldots,z_d)
	=
	\begin{cases}
		\displaystyle
		\sum_{\alpha\in\mathbb N_0^d\setminus\{0\}}
		c_\alpha z_1^{\alpha_1}\cdots z_d^{\alpha_d},
		& X\text{ commutative},\\[2mm]
		\displaystyle
		\sum_{r\ge1}\sum_{i_1,\ldots,i_r=1}^{d}
		c_{i_1,\ldots,i_r}z_{i_1}\cdots z_{i_r},
		& X\text{ noncommutative},
	\end{cases}
	\]
	where the coefficients belong to $\mathbb C$ and only finitely
	many are nonzero.
	Let $\mathcal M_d$ denote the corresponding monomial basis,
	and evaluate polynomials using $\ast$.
	
	A tuple $\mathbf u=(u_1,\ldots,u_d)\in X^d$
	\emph{freely generates} its algebra if
	\[
	P(\mathbf u)\ne0
	\qquad
	\bigl(P\in\mathscr P_d\setminus\{0\}\bigr).
	\]
	Equivalently, the evaluation homomorphism
	$P\mapsto P(\mathbf u)$ is injective.
	A sequence freely generates its algebra if every finite initial
	tuple does so. Freeness is understood in the commutative or
	noncommutative associative sense, according to $\mathscr P_d$.
	
	The following criterion adapts the fixed-monomial Baire argument
	of Bayart, Costa J\'unior and Papathanasiou%
	~\cite[Theorem~2.1 and Remark~2.2]{BCP2021}
	to finite-path actions.
	\begin{theorem}\label{thm:algebra-transitivity-criterion}
		Fix $d\ge1$.
		Suppose that every finite nonempty
		$A\subseteq\mathcal{M}_d$ admits $\beta_A\in A$ such that,
		for any nonempty open sets $U_1,\ldots,U_d,V\subseteq X$
		and any neighborhood $W$ of zero, there exist
		$\mathbf{u}\in U_1\times\cdots\times U_d$ and a nonempty
		finite path $\sigma$ satisfying
		\[
		T^\sigma\beta_A(\mathbf{u})\in V,
		\qquad
		T^\sigma\mu(\mathbf{u})\in W
		\quad(\mu\in A\setminus\{\beta_A\}).
		\]
		Then there exists a dense $G_\delta$ set
		$\mathcal{R}_d\subseteq X^d$ whose elements freely generate
		hypercyclic algebras for $\mathcal{T}$.
		
		If the hypothesis holds for every $d\ge1$, then
		$\mathcal{T}$ supports a dense hypercyclic algebra freely
		generated by countably infinitely many elements.
	\end{theorem}
	
	\begin{proof}
		Let $(V_j)_{j\ge1}$ be a countable base of nonempty open
		subsets of $X$.
		For a finite nonempty $A\subseteq\mathcal{M}_d$ and
		$s\ge1$, put
		\[
		K(A,s)=
		\left\{
		c=(c_\mu)_{\mu\in A}\in\mathbb{
			C}^{A}:
		c_{\beta_A}=1,\quad |c_\mu|\le s
		\right\},
		\]
		and write
		\[
		P_c(\mathbf{u})
		=
		\sum_{\mu\in A}c_\mu\mu(\mathbf{u}).
		\]
		Define
		\[
		\Omega(A,s,j)
		=
		\bigcup_{\sigma\ne\varnothing}
		\left\{
		\mathbf{u}\in X^d:
		T^\sigma P_c(\mathbf{u})\in V_j
		\text{ for every }c\in K(A,s)
		\right\},
		\]
		where the union ranges over all nonempty finite paths.
		
		For fixed $\sigma$, the map
		\[
		(\mathbf{u},c)\longmapsto T^\sigma P_c(\mathbf{u})
		\]
		is jointly continuous.
		Since $K(A,s)$ is compact, the set inside braces is open:
		a finite subcover of the coefficient neighborhoods gives
		a common neighborhood of any qualifying tuple.
		Hence $\Omega(A,s,j)$ is open.
		
		To prove density, fix nonempty open sets
		\(U_1,\ldots,U_d\subseteq X\) and choose \(v\in V_j\).
		Let \(r\) denote the number of elements of \(A\).
		Choose a balanced open neighborhood $D$ of zero such that
		\[
		v+\underbrace{D+\cdots+D}_{r\ \mathrm{summands}}
		\subseteq V_j.
		\]
		Apply the hypothesis with target $v+D$ and zero
		neighborhood $s^{-1}D$.
		It gives $\mathbf{u}\in U_1\times\cdots\times U_d$
		and a path $\sigma$ such that
		\[
		T^\sigma\beta_A(\mathbf{u})\in v+D,
		\qquad
		T^\sigma\mu(\mathbf{u})\in s^{-1}D
		\quad(\mu\ne\beta_A).
		\]
		For every $c\in K(A,s)$, linearity and balancedness yield
		\[
		\begin{aligned}
			T^\sigma P_c(\mathbf{u})
			=
			T^\sigma\beta_A(\mathbf{u})
			+
			\sum_{\mu\ne\beta_A}
			c_\mu T^\sigma\mu(\mathbf{u})\in
			v+\underbrace{D+\cdots+D}_{r\ \mathrm{summands}}
			\subseteq V_j.
		\end{aligned}
		\]
		Thus $\Omega(A,s,j)$ meets every nonempty open product
		and is therefore dense.
		This argument uses balanced neighborhoods only;
		local convexity is not required.
		
		The monomial set $\mathcal{M}_d$ is countable.
		By the Baire category theorem,
		\[
		\mathcal{R}_d
		=
		\bigcap_{\substack{
				A\subseteq\mathcal{M}_d\\
				A\ \mathrm{finite,\ nonempty}}}
		\ \bigcap_{s,j\ge1}\Omega(A,s,j)
		\]
		is a dense $G_\delta$ subset of $X^d$.
		
		Fix $\mathbf{u}\in\mathcal{R}_d$ and a nonzero polynomial
		\[
		P=\sum_{\mu\in A}a_\mu\mu,
		\qquad A=\operatorname{supp}P.
		\]
		All coefficients $a_\mu$, $\mu\in A$, are nonzero.
		In particular, $a_{\beta_A}\ne0$.
		Choose an integer $s\ge1$ such that
		\[
		s\ge
		\max_{\mu\in A}
		\left|\frac{a_\mu}{a_{\beta_A}}\right|,
		\]
		and set $Q=P/a_{\beta_A}$.
		Its coefficient vector belongs to $K(A,s)$.
		Since $\mathbf{u}\in\Omega(A,s,j)$ for every $j$,
		the orbit of $Q(\mathbf{u})$ meets every $V_j$.
		Therefore
		\[
		\overline{\operatorname{orb}
			(\mathcal{T},Q(\mathbf{u}))}=X.
		\]
		Every path operator is linear, so
		\[
		\overline{\operatorname{orb}
			(\mathcal{T},P(\mathbf{u}))}
		=
		a_{\beta_A}
		\overline{\operatorname{orb}
			(\mathcal{T},Q(\mathbf{u}))}
		=X.
		\]
		As $X\ne\{0\}$, this also implies $P(\mathbf{u})\ne0$.
		Thus $\mathbf{u}$ freely generates a hypercyclic algebra.
		
		Finally, suppose the hypothesis holds for every $d$.
		Let $\pi_d:X^{\mathbb{N}}\to X^d$ be the projection onto
		the first $d$ coordinates.
		The set
		\[
		\begin{aligned}
			\mathcal{G}
			=
			\bigcap_{d\ge1}\pi_d^{-1}(\mathcal{R}_d)\cap
			\bigcap_{j\ge1}\bigcup_{k\ge1}
			\left\{
			(x_h)_{h\ge1}:x_k\in V_j
			\right\}
		\end{aligned}
		\]
		is a dense $G_\delta$ subset of $X^{\mathbb{N}}$.
		Indeed, the projections are continuous and open, and
		each condition requiring some coordinate to belong to
		$V_j$ defines an open dense set.
		
		Every sequence in $\mathcal{G}$ has dense coordinate
		range, and each finite initial segment freely generates
		a hypercyclic algebra.
		Consequently, the entire sequence freely generates a
		dense hypercyclic algebra for $\mathcal{T}$.
	\end{proof}
	
	\begin{remark}\label{rem:fixed-monomial}
		The monomial $\beta_A$ is fixed before the open sets are chosen.
		In the applications below it is a monomial of uniquely minimal
		weighted degree, or a longest noncommutative word. This requirement
		allows the same criterion to treat noncommuting operator pairs.
	\end{remark}
	\section{Weighted shifts with coordinatewise multiplication}
	\label{sec:pointwise}
	
	We first fix notation used throughout the remaining sections.
	By a \emph{double Fr\'echet sequence space} we mean a linear
	subspace $X\subseteq\omega^2$ endowed with a Fr\'echet topology
	for which every coordinate projection
	$x\mapsto x_{i,j}$ is continuous.
	Write $e_{i,j}$ for the coordinate vectors and
	$F_N=\{0,\ldots,N\}^2$.
	The weighted backward shifts are
	\begin{equation}\label{eq:weighted-shifts}
		(T_1x)_{i,j}=w^{(1)}_{i+1,j}x_{i+1,j},
		\qquad
		(T_2x)_{i,j}=w^{(2)}_{i,j+1}x_{i,j+1}.
	\end{equation}
	The complex weights may vanish, and the operators need not commute.
	Their continuity on each space under consideration will be assumed.
	
	For a finite path $\sigma$, write its associated operator
	in block form as
	\[
	T^\sigma
	=
	T_1^{a_1}T_2^{b_1}\cdots T_1^{a_q}T_2^{b_q},
	\qquad
	q\ge1,\quad a_v,b_v\in\mathbb N_0,
	\]
	allowing zero exponents at the ends. Put
	\[
	A_0=B_0=0,\qquad
	A_v=\sum_{h=1}^{v}a_h,\qquad
	B_v=\sum_{h=1}^{v}b_h
	\qquad(1\le v\le q).
	\]
	The \emph{displacement} of $\sigma$ is the pair
	\[
	(m,n)=(A_q,B_q)
	=
	\left(
	\sum_{v=1}^{q}a_v,\,
	\sum_{v=1}^{q}b_v
	\right).
	\]
	Direct substitution gives
	\begin{equation}\label{eq:path-action}
		(T^\sigma x)_{i,j}
		=
		\Lambda_\sigma(i,j)x_{i+m,j+n},
	\end{equation}
	where
	\begin{equation}\label{eq:path-weight}
		\Lambda_\sigma(i,j)
		=
		\left(
		\prod_{v=1}^{q}\prod_{s=1}^{a_v}
		w^{(1)}_{i+A_{v-1}+s,\,j+B_{v-1}}
		\right)
		\left(
		\prod_{v=1}^{q}\prod_{t=1}^{b_v}
		w^{(2)}_{i+A_v,\,j+B_{v-1}+t}
		\right).
	\end{equation}
	Empty products equal $1$.
	All subsequent path conditions refer to these weights
	and displacements.
	
	\subsection{Sequence algebras with a continuous norm}
	
	Here multiplication is coordinatewise:
	$(x\cdot y)_{i,j}=x_{i,j}y_{i,j}$. A continuous norm on a Fr\'echet
	space need not generate its topology. The following characterization
	uses its existence together with a submultiplicative defining family.
	
	\begin{theorem}\label{thm:coordinatewise-shift-algebras}
		Let $X$ be a complex Fr\'echet sequence algebra on
		$\mathbb{N}_0^2$ with coordinatewise multiplication.
		Assume that its topology is generated by an increasing
		sequence of submultiplicative seminorms, that $X$ admits
		a continuous norm, and that the coordinate vectors form
		a Schauder basis in a fixed enumeration.
		Let $\mathcal T=(T_1,T_2)$ be continuous weighted shifts
		as in \eqref{eq:weighted-shifts}.
		Then the following conditions are equivalent:
		\begin{enumerate}
			\item[(i)]
			$\mathcal{T}$ supports a dense hypercyclic algebra freely
			generated by countably infinitely many elements,
			with freeness understood in the commutative sense.
			
			\item[(ii)]
			$\mathcal{T}$ supports a nonzero hypercyclic algebra.
			
			\item[(iii)]
			For every $r\ge1$, there exists $x\in X$ such that $x^r$
			is hypercyclic for $\mathcal{T}$.
			
			\item[(iv)]
			For every $r\ge1$ and $L\ge0$, there exists a sequence
			of finite paths $(\sigma_k)$, with displacements
			$(m_k,n_k)$ satisfying $m_k+n_k\to\infty$, such that
			\[
			|\Lambda_{\sigma_k}(i,j)|^{-1/r}
			e_{i+m_k,j+n_k}\longrightarrow0
			\qquad((i,j)\in F_L).
			\]
			The sequence may depend on $r$ and $L$.
			
			\item[(v)]
			There exists a sequence of finite paths $(\sigma_k)$,
			with displacements satisfying $m_k+n_k\to\infty$, such that
			\[
			|\Lambda_{\sigma_k}(i,j)|^{-\gamma}
			e_{i+m_k,j+n_k}\longrightarrow0
			\]
			for every $\gamma>0$ and every fixed
			$(i,j)\in\mathbb{N}_0^2$.
		\end{enumerate}
		All limits are taken in $X$.
		In \textup{(iv)} and \textup{(v)}, the relevant path weights
		are required to be eventually nonzero.
	\end{theorem}
	
	\begin{proof}
		Let $(p_q)$ be an increasing submultiplicative defining family.
		A continuous norm is bounded above by $C p_{q_0}$ for some
		$C>0$ and $q_0$, so $p_{q_0}$ is itself a norm. Discarding the
		preceding seminorms, we may assume that $p_1$ is a norm. Since $e_{i,j}^2=e_{i,j}$,
		submultiplicativity gives
		\[
		0<p_1(e_{i,j})\le p_1(e_{i,j})^2.
		\]
		Consequently,
		\begin{equation}\label{eq:coordinate-basis-lower-bound}
			p_q(e_{i,j})\ge p_1(e_{i,j})\ge1
			\qquad(q\ge1).
		\end{equation}
		
		The implication \textup{(i)}$\Rightarrow$\textup{(ii)}
		is immediate.
		For \textup{(ii)}$\Rightarrow$\textup{(iii)}, choose a nonzero
		element $x$ of a hypercyclic algebra.
		Coordinatewise multiplication has no nonzero nilpotent
		elements, so $x^r\ne0$ for every $r\ge1$.
		Thus every positive power of $x$ is hypercyclic.
		
		\smallskip
		\noindent
		\textup{(iii)} $\Rightarrow$ \textup{(iv)}.
		Fix $r\ge1$ and $L\ge0$, and choose $x\in X$ such that
		$x^r$ is hypercyclic for $\mathcal{T}$.
		We can choose paths $\sigma_k$ satisfying
		\[
		m_k+n_k\ge k,
		\qquad
		\max_{(i,j)\in F_L}
		\left|
		\Lambda_{\sigma_k}(i,j)
		x_{i+m_k,j+n_k}^{\,r}-1
		\right|
		<\frac{1}{k+1}.
		\]
		Indeed, there are only finitely many words of length
		less than $k$, and their orbit images can be excluded
		from the prescribed nonempty open neighborhood.
		
		For every fixed $(i,j)\in F_L$, the Schauder basis property
		gives
		\[
		x_{i+m_k,j+n_k}e_{i+m_k,j+n_k}\longrightarrow0.
		\]
		Hence, for each $q\ge1$,
		\[
		\begin{aligned}
			p_q\!\left(
			|\Lambda_{\sigma_k}(i,j)|^{-1/r}e_{i+m_k,j+n_k}
			\right)
			=
			\left|
			\Lambda_{\sigma_k}(i,j)x_{i+m_k,j+n_k}^{\,r}
			\right|^{-1/r}
			p_q\!\left(
			x_{i+m_k,j+n_k}e_{i+m_k,j+n_k}
			\right)
			\longrightarrow0.
		\end{aligned}
		\]
		This proves \textup{(iv)}.
		
		\smallskip
		\noindent
		\textup{(iv)} $\Rightarrow$ \textup{(v)}.
		Apply \textup{(iv)} with $r=k$ and $L=k$.
		Choose one path $\sigma_k$ from the resulting sequence
		such that
		\[
		m_k+n_k\ge k,
		\qquad
		\max_{(i,j)\in F_k}
		p_k\!\left(
		|\Lambda_{\sigma_k}(i,j)|^{-1/k}
		e_{i+m_k,j+n_k}
		\right)
		<\frac1k.
		\]
		By \eqref{eq:coordinate-basis-lower-bound},
		\[
		|\Lambda_{\sigma_k}(i,j)|>k^k
		\qquad((i,j)\in F_k).
		\]
		Fix $q\ge1$, $\gamma>0$, and $(i,j)\in\mathbb{N}_0^2$.
		For all sufficiently large $k$, we have
		$q,i,j\le k$ and $\gamma\ge1/k$. Therefore,
		\[
		\begin{aligned}
			p_q\!\left(
			|\Lambda_{\sigma_k}(i,j)|^{-\gamma}
			e_{i+m_k,j+n_k}
			\right)
			&\le
			p_k\!\left(
			|\Lambda_{\sigma_k}(i,j)|^{-1/k}
			e_{i+m_k,j+n_k}
			\right)\\
			&<\frac1k.
		\end{aligned}
		\]
		Thus the same path sequence satisfies \textup{(v)}.
		
		\smallskip
		\noindent
		\textup{(v)} $\Rightarrow$ \textup{(i)}.
		Fix a sequence as in \textup{(v)}.
		By \eqref{eq:coordinate-basis-lower-bound},
		\begin{equation}\label{eq:path-weights-diverge}
			|\Lambda_{\sigma_k}(i,j)|\longrightarrow\infty
		\end{equation}
		for every fixed $(i,j)$.
		
		We verify the monomial-separation hypothesis of the
		criterion in Theorem~\ref{thm:algebra-transitivity-criterion}.
		Fix $d\ge1$ and a finite nonempty set
		\[
		A\subseteq\mathbb{N}_0^d\setminus\{0\}.
		\]
		Choose $\kappa=(\kappa_1,\ldots,\kappa_d)\in\mathbb{N}^d$
		such that the integers
		\[
		D_\alpha=\kappa\cdot\alpha
		\qquad(\alpha\in A)
		\]
		are pairwise distinct.
		For example, sufficiently large successive powers of
		an integer may be used for the entries of $\kappa$.
		Let $\beta\in A$ be the unique minimizer and put
		$D=D_\beta$. Thus
		\[
		D_\alpha>D\qquad(\alpha\in A\setminus\{\beta\}).
		\]
		
		Let $U_1,\ldots,U_d,V$ be nonempty open subsets of $X$,
		and let $W$ be a neighborhood of zero.
		Choose finitely supported elements
		\[
		x_a\in U_a\quad(1\le a\le d),
		\qquad y\in V,
		\]
		with all supports contained in some $F_N$.
		For sufficiently large $k$, choose complex numbers
		$\eta_{i,j,k}$ satisfying
		\[
		\eta_{i,j,k}^{\,D}
		=
		\frac{y_{i,j}}{\Lambda_{\sigma_k}(i,j)}
		\qquad((i,j)\in F_N),
		\]
		and define
		\[
		u_{a,k}
		=
		x_a+
		\sum_{(i,j)\in F_N}
		\eta_{i,j,k}^{\,\kappa_a}
		e_{i+m_k,j+n_k},
		\qquad 1\le a\le d.
		\]
		For every $q\ge1$,
		\[
		\begin{aligned}
			p_q\!\left(
			\eta_{i,j,k}^{\,\kappa_a}
			e_{i+m_k,j+n_k}
			\right)
			=
			|y_{i,j}|^{\kappa_a/D}\times
			p_q\!\left(
			|\Lambda_{\sigma_k}(i,j)|^{-\kappa_a/D}
			e_{i+m_k,j+n_k}
			\right)
			\longrightarrow0.
		\end{aligned}
		\]
		Hence $u_{a,k}\to x_a$, so $u_{a,k}\in U_a$ eventually.
		
		When $m_k+n_k>2N$, the original and translated supports
		are disjoint, and $T^{\sigma_k}$ annihilates every vector
		supported in $F_N$.
		Coordinatewise multiplication and
		\eqref{eq:path-action} therefore give
		\begin{equation}\label{eq:monomial-image}
			T^{\sigma_k}(\mathbf{u}_k^\alpha)
			=
			\sum_{(i,j)\in F_N}
			\Lambda_{\sigma_k}(i,j)\eta_{i,j,k}^{\,D_\alpha}e_{i,j}
			\qquad(\alpha\in A).
		\end{equation}
		For $\alpha=\beta$, this is exactly $y$.
		For $\alpha\ne\beta$,
		\[
		\left|
		\Lambda_{\sigma_k}(i,j)\eta_{i,j,k}^{\,D_\alpha}
		\right|
		=
		|y_{i,j}|^{D_\alpha/D}
		|\Lambda_{\sigma_k}(i,j)|^{\,1-D_\alpha/D}
		\longrightarrow0
		\]
		by \eqref{eq:path-weights-diverge}.
		Since the sum in \eqref{eq:monomial-image} is finite,
		\[
		T^{\sigma_k}(\mathbf{u}_k^\beta)=y,
		\qquad
		T^{\sigma_k}(\mathbf{u}_k^\alpha)\longrightarrow0
		\quad(\alpha\in A\setminus\{\beta\}).
		\]
		Thus, for sufficiently large $k$,
		\[
		T^{\sigma_k}(\mathbf{u}_k^\beta)\in V,
		\qquad
		T^{\sigma_k}(\mathbf{u}_k^\alpha)\in W
		\quad(\alpha\ne\beta).
		\]
		The choice of $\beta$ depends only on $A$. Applying
		Theorem~\ref{thm:algebra-transitivity-criterion} for every $d\ge1$
		gives \textup{(i)}.
	\end{proof}
	Shu, Wang and Zhao~\cite{SWZ2017} characterized hypercyclicity
	for weighted double shifts, while B\`es and
	Papathanasiou~\cite{BP2020} established the existence of dense
	hypercyclic algebras for various convolution operators.
	For $1\le p<\infty$, write
	\[
	\mathcal L_p=\ell^p(\mathbb N_0^2),
	\qquad
	\|x\|_p=\Biggl(\sum_{i,j\ge0}|x_{i,j}|^p\Biggr)^{1/p}.
	\]
	For coordinatewise multiplication on $\mathcal L_p$
	and $c_0(\mathbb N_0^2)$, the following consequence of
	Theorem~\ref{thm:coordinatewise-shift-algebras} relates
	hypercyclicity to dense strong algebrability.
	
	\begin{corollary}\label{cor:lp-c0-hypercyclic-algebras}
		Let
		\[
		X=\mathcal L_p,\qquad 1\le p<\infty,
		\quad\text{or}\quad
		X=c_0(\mathbb N_0^2),
		\]
		with its usual norm and coordinatewise multiplication.
		Let $\mathcal T=(T_1,T_2)$ be a pair of bounded weighted
		backward shifts as in \eqref{eq:weighted-shifts}, not
		necessarily commuting. Then the following are equivalent:
		\begin{enumerate}
			\item[\textup{(i)}]
			$\mathcal T$ is hypercyclic.
			
			\item[\textup{(ii)}]
			$\mathcal T$ supports a nonzero hypercyclic algebra.
			
			\item[\textup{(iii)}]
			$\mathcal T$ supports a dense hypercyclic algebra freely
			generated by countably infinitely many elements.
			
			\item[\textup{(iv)}]
			There exists a sequence of finite paths $(\sigma_k)$,
			with displacements $(m_k,n_k)$ satisfying
			$m_k+n_k\to\infty$, such that
			\[
			\min_{(i,j)\in F_N}
			|\Lambda_{\sigma_k}(i,j)|
			\longrightarrow\infty
			\qquad\text{for every }N\ge0.
			\]
		\end{enumerate}
	\end{corollary}
	
	\begin{proof}
		Both spaces are Banach sequence algebras with
		submultiplicative norms and Schauder coordinate bases,
		and $\|e_{i,j}\|=1$. Hence
		Theorem~\ref{thm:coordinatewise-shift-algebras} applies.
		Clearly, \textup{(iii)} implies \textup{(ii)}, which
		implies \textup{(i)}.
		
		The implication \textup{(i)} $\Rightarrow$ \textup{(iv)}
		follows from the necessity argument in
		\cite[Theorem~4.3]{SWZ2017}.
		Although stated there for positive weights on $\mathcal L_p$,
		the same coordinate estimate applies to $c_0$ and to
		complex weights, possibly zero, by using absolute values
		and arbitrarily small nonzero scalar multiples of a
		hypercyclic vector. Boundedness of the shifts ensures
		that the resulting displacements satisfy
		$m_k+n_k\to\infty$.
		
		Finally, assume \textup{(iv)}. For every $\gamma>0$
		and every fixed $(i,j)$,
		\[
		\left\|
		|\Lambda_{\sigma_k}(i,j)|^{-\gamma}
		e_{i+m_k,j+n_k}
		\right\|
		=
		|\Lambda_{\sigma_k}(i,j)|^{-\gamma}
		\longrightarrow0.
		\]
		Thus condition \textup{(v)} of
		Theorem~\ref{thm:coordinatewise-shift-algebras}
		holds, yielding \textup{(iii)}.
	\end{proof}
	Bourdon and Shapiro~\cite{BourdonShapiro2000} studied
	operators commuting with the backward Taylor shift
	\[
	(Bf)(z)=\frac{f(z)-f(0)}{z}
	\]
	on $A^2(\mathbb D)$.
	Let $\mathbb D^2=\{(z,w)\in\mathbb C^2:|z|<1,\ |w|<1\}$.
	The Bergman space $A^2(\mathbb D^2)$ consists of holomorphic
	functions $f(z,w)=\sum_{m,n\ge0}a_{m,n}z^mw^n$ such that
	\[
	\|f\|_{A^2}^2
	=
	\int_{\mathbb D^2}|f(z,w)|^2\,dA(z)\,dA(w)
	=
	\sum_{m,n\ge0}\frac{|a_{m,n}|^2}{(m+1)(n+1)}
	<\infty,
	\]
	where $dA(\zeta)=\pi^{-1}\,dx\,dy$, $\zeta=x+iy$.
	The coefficient identity follows from the orthogonality
	of the monomials and
	$\int_{\mathbb D}|\zeta|^{2m}\,dA(\zeta)=(m+1)^{-1}$.
	The coordinate Taylor backward shifts are
	\[
	(T_1f)(z,w)=\frac{f(z,w)-f(0,w)}{z},
	\qquad
	(T_2f)(z,w)=\frac{f(z,w)-f(z,0)}{w},
	\]
	with holomorphic extensions across the coordinate axes.
	
	\begin{example}\label{ex:bergman-noncommuting-pair}
		Let $\mathcal H=A^2(\mathbb D^2)$ carry normalized area
		measure. For $f(z,w)=\sum_{m,n\ge0}a_{m,n}z^mw^n$ and $g(z,w)=\sum_{m,n\ge0}b_{m,n}z^mw^n$, define
		\[
		(f\diamond g)(z,w)
		=
		\sum_{m,n\ge0}
		\frac{a_{m,n}b_{m,n}}{\sqrt{(m+1)(n+1)}}z^mw^n.
		\]
		The unitary map
		\[
		U:\mathcal H\longrightarrow\mathcal L_2,
		\qquad
		Uf=
		\left(
		\frac{a_{m,n}}{\sqrt{(m+1)(n+1)}}
		\right)_{m,n\ge0},
		\]
		satisfies
		\[
		U(f\diamond g)=(Uf)\cdot(Ug).
		\]
		Thus $(\mathcal H,\diamond)$ is a commutative Banach algebra,
		with
		\[
		\|f\diamond g\|_{\mathcal H}
		\le \|f\|_{\mathcal H}\|g\|_{\mathcal H}.
		\]
		
		Consider the coordinate Taylor backward shifts
		\[
		\begin{aligned}
			(T_1f)(z,w)
			=\sum_{m,n\ge0}a_{m+1,n}z^mw^n,\quad
			(T_2f)(z,w)
			=\sum_{m,n\ge0}a_{m,n+1}z^mw^n.
		\end{aligned}
		\]
		These are commuting coordinate versions of the Bergman
		backward shift studied by Bourdon and
		Shapiro~\cite{BourdonShapiro2000}.
		Writing $\widetilde T_j=UT_jU^{-1}$, we obtain
		\[
		\begin{aligned}
			(\widetilde T_1x)_{i,j}
			=\sqrt{\frac{i+2}{i+1}}\,x_{i+1,j},\quad
			(\widetilde T_2x)_{i,j}
			=\sqrt{\frac{j+2}{j+1}}\,x_{i,j+1}.
		\end{aligned}
		\]
		In particular, $\|T_1\|,\|T_2\|\le\sqrt2$.
		
		For the paths
		\[
		\widetilde T^{\sigma_k}=\widetilde T_2^k,
		\qquad (m_k,n_k)=(0,k),
		\]
		the weight products are
		\[
		\Lambda_{\sigma_k}(i,j)
		=
		\prod_{s=1}^{k}\sqrt{\frac{j+s+1}{j+s}}
		=
		\sqrt{\frac{j+k+1}{j+1}}.
		\]
		Hence, for every $N\ge0$,
		\[
		\min_{(i,j)\in F_N}|\Lambda_{\sigma_k}(i,j)|
		=
		\sqrt{1+\frac{k}{N+1}}
		\longrightarrow\infty.
		\]
		By Corollary~\ref{cor:lp-c0-hypercyclic-algebras}, the
		conjugate pair supports a dense hypercyclic algebra freely
		generated by countably infinitely many elements.
		Transferring this algebra through $U^{-1}$ gives the same
		conclusion for $\mathcal T=(T_1,T_2)$ on
		$(\mathcal H,\diamond)$.
	\end{example}
\begin{remark}
	The choice of product in
	Example~\ref{ex:bergman-noncommuting-pair}
	is essential for the hypercyclic-algebra conclusion.
	Consider instead
	\[
	(f\circ g)(z,w)
	=
	\sum_{m,n\ge0}
	\frac{a_{m,n}b_{m,n}}{(m+1)(n+1)}z^mw^n,
	\]
	where
	$f(z,w)=\sum_{m,n\ge0}a_{m,n}z^mw^n$ and
	$g(z,w)=\sum_{m,n\ge0}b_{m,n}z^mw^n$.
	This product is associative and commutative and satisfies
	\[
	\|f\circ g\|_{\mathcal H}
	\le
	\|f\|_{\mathcal H}\|g\|_{\mathcal H}.
	\]
	Thus $(\mathcal H,\circ)$ is a commutative Banach algebra
	with its usual Bergman norm.
	We show that the same pair of coordinate Taylor backward
	shifts supports no nonzero hypercyclic algebra for
	this product.
	
	Since $T_1$ and $T_2$ commute, every finite composition
	has the form $T_1^kT_2^l$, with $k,l\ge0$.
	For every $f\in\mathcal H$, we have
	\[
	(f\circ f)(z,w)
	=
	\sum_{m,n\ge0}
	\frac{a_{m,n}^{\,2}}{(m+1)(n+1)}z^mw^n.
	\]
	Consequently, for all $k,l\ge0$,
	\[
	\begin{aligned}
		\left|
		\bigl(T_1^kT_2^l(f\circ f)\bigr)(0,0)
		\right|
		&=
		\frac{|a_{k,l}|^2}{(k+1)(l+1)}\\
		&\le
		\sum_{m,n\ge0}
		\frac{|a_{m,n}|^2}{(m+1)(n+1)}
		=\|f\|_{\mathcal H}^2.
	\end{aligned}
	\]
	Thus the image of the orbit of $f\circ f$ under
	evaluation at $(0,0)$ is bounded.
	This evaluation map is continuous and surjective
	from $\mathcal H$ onto $\mathbb C$, so a dense orbit
	would have a dense image in $\mathbb C$.
	Hence $f\circ f$ is not hypercyclic for $\mathcal T$.
	
	Finally, $f\ne0$ implies $f\circ f\ne0$.
	Every nonzero subalgebra therefore contains a nonzero
	square that is not hypercyclic.
	It follows that $\mathcal T$ supports no nonzero
	hypercyclic algebra on $(\mathcal H,\circ)$.
\end{remark}
	
	\subsection{The product topology on
		\texorpdfstring{$\omega^2$}{omega squared}}
	On $\omega^2$, the product topology is generated by the seminorms
	(see \cite[Section~3]{Wang2026})
	\[
	p_N(x)=\max_{(i,j)\in F_N}|x_{i,j}|,\qquad N\ge0.
	\]
	Thus convergence means convergence at each coordinate.
	Every continuous seminorm is bounded above by a multiple
	of some $p_N$ and therefore vanishes on arrays supported
	outside $F_N$. Consequently, $\omega^2$ admits no continuous
	norm, so Theorem~\ref{thm:coordinatewise-shift-algebras}
	does not apply.
	The following result adapts the compact-coefficient argument
	of \cite[Theorem~4.8]{BCP2021} from single weighted shifts
	to operator pairs on $\omega^2$, allowing vanishing weights.

	\begin{theorem}\label{prop:omega2-pointwise-algebras}
		Let $X=\omega^2$ carry the product topology and
		coordinatewise multiplication. For a pair
		$\mathcal T=(T_1,T_2)$ of weighted backward shifts
		as above, the following are equivalent:
		\begin{enumerate}
			\item[\textup{(i)}]
			$\mathcal T$ supports a nonzero hypercyclic algebra.
			
			\item[\textup{(ii)}]
			$\mathcal T$ supports a dense hypercyclic algebra.
			
			\item[\textup{(iii)}]
			At least one of $T_1$ and $T_2$ is hypercyclic.
			
			\item[\textup{(iv)}]
			At least one effective weight array is nowhere zero:
			\[
			w^{(1)}_{i+1,j}\ne0\quad\text{for all }i,j\ge0,
			\qquad\text{or}\qquad
			w^{(2)}_{i,j+1}\ne0\quad\text{for all }i,j\ge0.
			\]
		\end{enumerate}
		Whenever these conditions hold, the dense hypercyclic
		algebra can be chosen to be singly generated.
	\end{theorem}
	
	\begin{proof}
		Clearly, \textup{(ii)} implies \textup{(i)}.
		We first show that each of \textup{(i)} and
		\textup{(iii)} implies \textup{(iv)}.
		
		Either assumption implies that $\mathcal T$ is
		hypercyclic. Suppose that \textup{(iv)} fails.
		Choose $(i_1,j_1),(i_2,j_2)\in\mathbb N_0^2$ such that
		\[
		w^{(1)}_{i_1+1,j_1}=0,
		\qquad
		w^{(2)}_{i_2,j_2+1}=0,
		\]
		and define the proper closed coordinate hyperplanes
		\[
		H_1=\{z\in X:z_{i_1,j_1}=0\},
		\qquad
		H_2=\{z\in X:z_{i_2,j_2}=0\}.
		\]
		Then $T_1(X)\subseteq H_1$ and $T_2(X)\subseteq H_2$.
		Every nontrivial composition has $T_1$ or $T_2$ as its
		leftmost factor, so
		\[
		\operatorname{orb}(\mathcal T,x)
		\subseteq \{x\}\cup H_1\cup H_2
		\qquad(x\in X).
		\]
		The right-hand side is a proper closed subset of $X$.
		Hence $\mathcal T$ is not hypercyclic, a contradiction.
		This proves both implications.
		
		We now prove that \textup{(iv)} implies
		\textup{(ii)} and \textup{(iii)}.
		Assume the first array is nowhere zero and write $T=T_1$.
		For integers $1\le p\le q$ and $s\ge1$, put
		\[
		K_{p,q,s}
		=
		\left\{
		P(t)=\sum_{r=p}^q a_rt^r:
		|a_r|\le s,\quad |a_p|,|a_q|\ge s^{-1}
		\right\}.
		\]
		These compact families cover all nonzero polynomials
		with zero constant term. For such a family $K_{p,q,s}$
		and a nonempty open set $V$, define
		\[
		\mathcal E(K_{p,q,s},V)
		=
		\{x\in X:\ \forall P\in K_{p,q,s}\ \exists h\ge1,
		\ T^hP(x)\in V\}.
		\]
		Joint continuity of $(x,P)\mapsto T^hP(x)$ and a finite
		subcover of $K_{p,q,s}$ show that
		$\mathcal E(K_{p,q,s},V)$ is open.
		
		To prove density, fix a nonempty open set $U$.
		Choose finitely supported $f\in U$ and $y\in V$,
		and then $N$ and $\varepsilon>0$ such that their
		supports lie in $F_N$ and
		\[
		z|_{F_N}=f|_{F_N}\ \Longrightarrow\ z\in U,
		\qquad
		\max_{(i,j)\in F_N}|z_{i,j}-y_{i,j}|<\varepsilon
		\ \Longrightarrow\ z\in V.
		\]
		Choose $h_1>N$ and $h_{k+1}>h_k+N$.
		The blocks $F_N+(h_k,0)$ are disjoint, and the
		nonzero numbers
		\[
		\lambda_k(i,j)
		=
		\prod_{r=1}^{h_k}w^{(1)}_{i+r,j}
		\qquad((i,j)\in F_N)
		\]
		have, after passage to a common subsequence, limits
		$c_{i,j}\in\widehat{\mathbb C}$, where
		$\widehat{\mathbb C}=\mathbb C\cup\{\infty\}$
		denotes the extended complex plane.
		
		For $P(t)=\sum_{r=p}^q a_rt^r\in K_{p,q,s}$, choose
		\[
		\zeta_{k,i,j}(P)
		=
		\begin{cases}
			\bigl(y_{i,j}/(a_p\lambda_k(i,j))\bigr)^{1/p},
			&c_{i,j}=\infty,\\
			\bigl(y_{i,j}/(a_q\lambda_k(i,j))\bigr)^{1/q},
			&c_{i,j}=0,\\
			\zeta_{i,j},
			&0<|c_{i,j}|<\infty,
		\end{cases}
		\]
		where $P(\zeta_{i,j})=y_{i,j}/c_{i,j}$ in the last case.
		Roots can be chosen arbitrarily, and zero is chosen
		when $y_{i,j}=0$.
		If $c_{i,j}=\infty$, then
		\[
		\lambda_k(i,j)a_p\zeta_{k,i,j}(P)^p=y_{i,j},
		\]
		while, for $p<r\le q$,
		\[
		|\lambda_k(i,j)\zeta_{k,i,j}(P)^r|
		=
		\left|\frac{y_{i,j}}{a_p}\right|^{r/p}
		|\lambda_k(i,j)|^{1-r/p}
		\longrightarrow0.
		\]
		The case $c_{i,j}=0$ is analogous, using the term of
		degree $q$, and the finite nonzero case follows by
		continuity. Since $F_N$ is finite, these observations give
		\[
		\lambda_k(i,j)P(\zeta_{k,i,j}(P))
		\longrightarrow y_{i,j},
		\qquad
		C_P:=
		\max\biggl\{
		1,
		\sup_{\substack{k\ge1,\ (i,j)\in F_N\\p\le r\le q}}
		|\lambda_k(i,j)\zeta_{k,i,j}(P)^r|
		\biggr\}<\infty.
		\]
		
		Choose $k(P)$ so that the first limit has error less
		than $\varepsilon/2$ on $F_N$ for $k\ge k(P)$.
		If $Q(t)=\sum_{r=p}^q b_rt^r$ satisfies
		\[
		\max_{p\le r\le q}|b_r-a_r|
		<
		\frac{\varepsilon}{2(q-p+1)C_P},
		\]
		then, for all such $k$,
		\[
		\max_{(i,j)\in F_N}
		|\lambda_k(i,j)Q(\zeta_{k,i,j}(P))-y_{i,j}|
		<\varepsilon.
		\]
		Cover $K_{p,q,s}$ by finitely many of these coefficient
		neighborhoods, with centers $P_1,\ldots,P_t$, and choose
		distinct indices $k_\nu\ge k(P_\nu)$. Set
		\[
		u=f+\sum_{\nu=1}^t\sum_{(i,j)\in F_N}
		\zeta_{k_\nu,i,j}(P_\nu)e_{i+h_{k_\nu},j}.
		\]
		Then $u\in U$. For every $Q\in K_{p,q,s}$, a suitable
		center $P_\nu$ and the disjointness of the blocks give
		\[
		(T^{h_{k_\nu}}Q(u))_{i,j}
		=
		\lambda_{k_\nu}(i,j)
		Q(\zeta_{k_\nu,i,j}(P_\nu))
		\qquad((i,j)\in F_N).
		\]
		Hence $T^{h_{k_\nu}}Q(u)\in V$, proving density of
		$\mathcal E(K_{p,q,s},V)$.
		
		Intersecting these open dense sets over $p,q,s$ and
		a countable base of targets gives a dense $G_\delta$
		set $\mathcal R$ such that
		\[
		P(x)\in HC(T)
		\qquad
		\bigl(x\in\mathcal R,\
		P\in\mathbb C[t]\setminus\{0\},\ P(0)=0\bigr).
		\]
		
		For each $N\ge0$, let $\mathcal D_N$ be the set of
		arrays whose coordinates indexed by $F_N$ are nonzero
		and pairwise distinct. This set is open, since its
		defining conditions require finitely many continuous
		coordinate functions and their differences to be
		nonzero.
		
		To prove density of $\mathcal D_N$, fix a nonempty basic
		open set $U$ and choose $M\ge N$ such that all coordinates
		restricted by $U$ lie in $F_M$.
		For $(i,j)\in F_M$, choose $\xi_{i,j}$ successively,
		satisfying these restrictions and avoiding zero and
		the previously chosen values. Each choice is possible
		because a nonempty open subset of $\mathbb C$ remains
		nonempty after removing finitely many points.
		Define $u\in X$ by
		\[
		u_{i,j}=
		\begin{cases}
			\xi_{i,j},&(i,j)\in F_M,\\
			0,&(i,j)\notin F_M.
		\end{cases}
		\]
		Then $u\in U\cap\mathcal D_M\subseteq U\cap\mathcal D_N$,
		proving density.
		
		By the Baire category theorem, choose
		\[
		x\in\mathcal R\cap\bigcap_{N\ge0}\mathcal D_N.
		\]
		All coordinates of $x$ are nonzero and pairwise distinct.
		Put
		\[
		\mathcal A
		=
		\{P(x):P\in\mathbb C[t],\ P(0)=0\}.
		\]
		This is the nonunital algebra generated by $x$.
		Taking $P(t)=t$ gives $x\in HC(T)$, proving
		\textup{(iii)} and ensuring that $\mathcal A$ is nonzero.
		Every nonzero element of $\mathcal A$ is hypercyclic
		for $T$, and hence for $\mathcal T$.
		
		To prove density of $\mathcal A$, let $U\subseteq X$
		be nonempty and open. Choose $y\in U$ and $N\ge0$
		such that
		\[
		z|_{F_N}=y|_{F_N}\ \Longrightarrow\ z\in U.
		\]
		The numbers $0$ and $x_{i,j}$, $(i,j)\in F_N$,
		are distinct. Hence the interpolation polynomial
		\[
		P(t)
		=
		\sum_{(i,j)\in F_N}
		y_{i,j}\frac{t}{x_{i,j}}
		\prod_{\substack{(a,b)\in F_N\\(a,b)\ne(i,j)}}
		\frac{t-x_{a,b}}{x_{i,j}-x_{a,b}}
		\]
		is well defined and satisfies
		\[
		P(0)=0,\qquad
		P(x_{i,j})=y_{i,j}\quad((i,j)\in F_N).
		\]
		Thus $P(x)\in\mathcal A\cap U$.
		Since $U$ was arbitrary, $\mathcal A$ is dense in $X$.
		This proves \textup{(ii)} and the additional assertion.
		
		The second alternative in \textup{(iv)} follows by
		interchanging the two indices.
	\end{proof}
	\section{Hypercyclic algebras for weighted double shifts with the triangular product}
	\label{sec:triangular-product}
	
	We first identify the topology naturally associated with the triangular
	product. The underlying space will have the same product topology as in
	the preceding subsection.
	
	\begin{lemma}\label{lem:triangular-frechet-algebra}
		Let $X=\omega^2$ and define
		\[
		(x\ast y)_{i,j}=\sum_{t=0}^{j}x_{i,t}y_{i+t,j-t},
		\qquad
		q_N(x)=\sum_{i+j\le N}|x_{i,j}| \quad(N\ge0).
		\]
		Then $(X,\ast)$ is a separable noncommutative Fr\'echet algebra.
		The increasing seminorms $(q_N)$ are submultiplicative and generate
		the product topology.
	\end{lemma}
	
	\begin{proof}
		For $N\ge0$, define the $(N+1)\times(N+1)$ upper triangular matrix
		\[
		\pi_N(x)=
		\begin{pmatrix}
			x_{0,0} & x_{0,1} & \cdots & x_{0,N}\\
			0       & x_{1,0} & \cdots & x_{1,N-1}\\
			\vdots  & \ddots  & \ddots & \vdots\\
			0       & \cdots  & 0      & x_{N,0}
		\end{pmatrix}.
		\]
		For $0\le a\le b\le N$, the multiplication formula gives
		\[
		\begin{aligned}
			(\pi_N(x)\pi_N(y))_{a,b}
			=\sum_{c=a}^{b}x_{a,c-a}y_{c,b-c}=(x\ast y)_{a,b-a}
			=(\pi_N(x\ast y))_{a,b}.
		\end{aligned}
		\]
		Thus $\pi_N(x\ast y)=\pi_N(x)\pi_N(y)$.
		Since the maps $\pi_N$ separate points, $\ast$ is associative.
		Writing $\|\cdot\|_\Sigma$ for the sum of the absolute values of all
		matrix entries, we obtain
		\[
		q_N(x\ast y)
		=\|\pi_N(x)\pi_N(y)\|_\Sigma
		\le q_N(x)q_N(y).
		\]
		For the seminorms $(p_N)$ introduced above,
		\[
		p_N(x)\le q_{2N}(x),
		\qquad
		q_N(x)\le (N+1)^2p_N(x).
		\]
		Thus $(q_N)$ generates the same product topology.
		As a countable product of copies of $\mathbb C$,
		$X$ is a separable Fr\'echet space in this topology.
		
		For $x,y,h,k\in X$, submultiplicativity gives
		\[
		q_N\bigl((x+h)\ast(y+k)-x\ast y\bigr)
		\le
		q_N(x)q_N(k)+q_N(h)q_N(y)+q_N(h)q_N(k).
		\]
		The right-hand side tends to zero as $h,k\to0$
		for every $N$, proving joint continuity of multiplication.
		Finally,
		\[
		e_{0,0}\ast e_{0,1}=e_{0,1},\qquad
		e_{0,1}\ast e_{0,0}=0,
		\]
		so the algebra is noncommutative.
	\end{proof}
	
	Under the identification $(i,j)\leftrightarrow[i,i+j]$,
	the product $\ast$ corresponds to the incidence convolution
	of the chain $\mathbb N_0$.
	For $X=\omega^2$ or $X=\mathcal L_p$ with $1\le p\le2$,
	define
	\[
	J_s=\{x\in X:x_{i,j}=0\text{ whenever }j<s\},
	\qquad s\ge1.
	\]
	The multiplication formula gives
	$J_s\ast J_t\subseteq J_{s+t}$ for $s,t\ge1$.
	It also shows that every commutator
	$x\ast y-y\ast x$, with $x,y\in X$, belongs to $J_1$.
	
	\subsection{A path characterization on
		\texorpdfstring{$\omega^2$}{omega squared}}
	
	We first consider $(\omega^2,\ast)$ with the product
	topology generated by $(q_N)$.
	The following theorem relates hypercyclic powers of
	zero-diagonal elements to dense free hypercyclic algebras.
	Both properties are characterized by nonvanishing path
	weights on finite rectangles and arbitrarily large
	second-coordinate displacements.

	\begin{theorem}\label{thm:triangular-shift-algebras}
		Let $\mathcal T=(T_1,T_2)$ be a pair of weighted backward shifts on
		$(\omega^2,\ast)$, with possibly vanishing complex weights.
		The following conditions are equivalent:
		\begin{enumerate}
			\item[\textup{(i)}]
			For every $r\ge1$, there exists $x\in J_1$ such that $x^{\ast r}$
			is hypercyclic for $\mathcal T$.
			
			\item[\textup{(ii)}]
			The pair $\mathcal T$ supports a dense hypercyclic algebra freely
			generated, as a nonunital associative algebra, by countably
			infinitely many elements.
			
			\item[\textup{(iii)}]
			For every $N\ge0$ and $L\ge1$, there exists a finite path $\sigma$
			with displacement $(m,n)$ such that $n\ge L$ and
			\begin{equation}\label{eq:triangular-path-condition}
				\Lambda_\sigma(i,j)\ne0\qquad((i,j)\in F_N).
			\end{equation}
			
			\item[\textup{(iv)}]
			The pair $\mathcal T$ supports a dense hypercyclic algebra.
		\end{enumerate}
	\end{theorem}
	
	\begin{proof}
		\textup{(i)} $\Rightarrow$ \textup{(iii)}.
		Fix $N,L$ and choose $x\in J_1$ such that
		$z=x^{\ast L}$ is hypercyclic.
		Choose a path $\sigma$, with displacement $(m,n)$, satisfying
		\[
		q_{2N}\biggl(
		T^\sigma z-\sum_{(i,j)\in F_N}e_{i,j}
		\biggr)<\tfrac12.
		\]
		Its $(0,0)$ coordinate satisfies
		\[
		(T^\sigma z)_{0,0}
		=\Lambda_\sigma(0,0)z_{m,n}\ne0.
		\]
		Since $z\in J_L$ implies $z_{m,n}=0$ whenever $n<L$,
		the nonvanishing of $z_{m,n}$ forces $n\ge L$.
		Every output coordinate on $F_N$ is nonzero, so
		\eqref{eq:path-action} gives
		\eqref{eq:triangular-path-condition}.
		
		\textup{(iii)} $\Rightarrow$ \textup{(ii)}.
		We verify Theorem~\ref{thm:algebra-transitivity-criterion}.
		Fix $d\ge1$ and a finite nonempty set $A$ of nonempty words, and
		choose a longest word $\beta=z_{\beta_1}\cdots z_{\beta_\ell}\in A$.
		This choice depends only on $A$.
		Given nonempty open sets $U_1,\ldots,U_d,V$ and a zero neighborhood
		$W$, choose finitely supported $x_a\in U_a$ and $y\in V$.
		Choose $N$ containing their supports and $\varepsilon>0$ such that,
		for $B=\{z:q_N(z)<\varepsilon\}$,
		\[
		x_a+B\subseteq U_a,
		\qquad y+B\subseteq V,
		\qquad B\subseteq W.
		\]
		By \textup{(iii)}, choose a path $\sigma$ with
		\[
		n>2N+(N+1)^2(\ell-1),\qquad
		\Lambda_\sigma(i,j)\ne0\quad((i,j)\in F_N).
		\]
		For $(i,j)\in F_N$, set
		\[
		p_0(i,j)=m+i,\qquad p_\ell(i,j)=m+n+i+j,
		\]
		and, for $1\le t<\ell$, set
		\[
		p_t(i,j)=m+2N+\bigl((N+1)i+j\bigr)(\ell-1)+t.
		\]
		Each sequence is strictly increasing. All its internal indices are
		distinct from those of the other sequences and lie in $(m+2N,m+n)$.
		There are no internal indices when $\ell=1$.
		
		Writing $E_{a,b}=e_{a,b-a}$ for $a\le b$, define
		\[
		\begin{aligned}
			u_a=x_a+\sum_{(i,j)\in F_N}\Bigg[
			\sum_{t=1}^{\ell-1}\delta_{a,\beta_t}
			E_{p_{t-1}(i,j),p_t(i,j)}+\delta_{a,\beta_\ell}
			\frac{y_{i,j}}{\Lambda_\sigma(i,j)}
			E_{p_{\ell-1}(i,j),p_\ell(i,j)}\Bigg].
		\end{aligned}
		\]
		Every added entry has terminal index greater than $2N$, so
		$q_N(u_a-x_a)=0$ and $u_a\in U_a$.
		
		We record the coordinate calculation on the whole space for later use.
		Let $w=z_{c_1}\cdots z_{c_r}\in A$, so $r\le\ell$.
		For any $(i,j)\in\mathbb N_0^2$,
		\begin{equation}\label{eq:word-coordinate-expansion}
			\bigl(w(\mathbf u)\bigr)_{m+i,n+j}
			=\sum_{\substack{q_0=m+i,\ q_r=m+n+i+j\\
					q_0\le\cdots\le q_r}}
			\prod_{t=1}^{r}(u_{c_t})_{q_{t-1},q_t-q_{t-1}}.
		\end{equation}
		All entries in the construction have terminal index at most
		$m+n+2N$. A nonzero summand therefore satisfies
		\[
		m\le q_0\le m+2N<m+n\le q_r\le m+n+2N.
		\]
		Set $t_0=\min\{t:q_t>m+2N\}$.
		Original entries have terminal index at most $2N$, and every added
		entry other than a first entry starts above $m+2N$. Hence, for some
		$(i',j')\in F_N$,
		\[
		(q_{t_0-1},q_{t_0})=(p_0(i',j'),p_1(i',j')).
		\]
		At every internal index $p_s=p_s(i',j')$, the definition gives
		\[
		(u_a)_{p_s,b-p_s}=0
		\quad\text{if }b\ge p_s,\ b\ne p_{s+1}
		\qquad(1\le s<\ell).
		\]
		Thus successive nonzero factors force
		\[
		q_{t_0+s-1}=p_s(i',j')
		\quad\bigl(1\le s\le\min\{\ell,r-t_0+1\}\bigr).
		\]
		If $r-t_0+1<\ell$, this gives $q_r<m+n$, a contradiction.
		Consequently,
		\[
		\ell\le r-t_0+1\le r\le\ell,
		\]
		so $r=\ell$ and $t_0=1$. Comparing endpoints gives
		$(i,j)=(i',j')\in F_N$.
		Substitution into \eqref{eq:word-coordinate-expansion} yields
		\[
		\bigl(w(\mathbf u)\bigr)_{m+i,n+j}
		=\begin{cases}
			\displaystyle\frac{y_{i,j}}{\Lambda_\sigma(i,j)}
			\prod_{t=1}^{\ell}\delta_{c_t,\beta_t},
			&r=\ell,\ (i,j)\in F_N,\\[2mm]
			0,&\text{otherwise}.
		\end{cases}
		\]
		Together with \eqref{eq:path-action}, this proves the global identities
		\begin{equation}\label{eq:triangular-word-separation}
			T^\sigma\beta(\mathbf u)=y,\qquad
			T^\sigma w(\mathbf u)=0\quad(w\in A\setminus\{\beta\}).
		\end{equation}
		The Baire criterion applies for every $d$ and proves \textup{(ii)}.
		
		\textup{(ii)} $\Rightarrow$ \textup{(iv)} is immediate.
		
		\textup{(iv)} $\Rightarrow$ \textup{(i)}.
		Let $\mathcal A$ be a dense hypercyclic algebra and fix $r\ge1$.
		Put
		\[
		D=\sum_{s=0}^{r}s e_{s,0},
		\qquad
		S=\sum_{s=0}^{r-1}e_{s,1}.
		\]
		The multiplication formula gives
		\[
		D\ast S-S\ast D
		=\sum_{s=0}^{r-1}\bigl(s-(s+1)\bigr)e_{s,1}
		=-S,
		\qquad
		S^{\ast r}=e_{0,r}.
		\]
		Hence
		\[
		(D\ast S-S\ast D)^{\ast r}
		=(-1)^r e_{0,r}\ne0.
		\]
		By density of $\mathcal A$ and continuity of multiplication
		and of the coordinate projections, choose $a,b\in\mathcal A$ such that
		\[
		\left|
		\bigl((a\ast b-b\ast a)^{\ast r}\bigr)_{0,r}
		-(-1)^r
		\right|<\tfrac12.
		\]
		Set $h=a\ast b-b\ast a$.
		Since
		\[
		h_{i,0}=a_{i,0}b_{i,0}-b_{i,0}a_{i,0}=0
		\qquad(i\ge0),
		\]
		we have $h\in\mathcal A\cap J_1$.
		The preceding inequality gives $h^{\ast r}\ne0$.
		Thus $h^{\ast r}$ is hypercyclic, since it belongs to
		$\mathcal A\setminus\{0\}$.
		As $r\ge1$ was arbitrary, \textup{(i)} follows.
	\end{proof}

\subsection{Banach sequence algebras and dependence on multiplication}
	
The triangular product makes $\mathcal L_p$, with its
usual norm, a Banach algebra for $1\le p\le2$.
The preceding matrix-unit construction then gives the
following norm characterization.

\begin{proposition}\label{thm:l1-triangular-algebras}
	Let $1\le p\le2$, and let $\mathcal T=(T_1,T_2)$ be
	a pair of bounded weighted backward shifts on
	$(\mathcal L_p,\ast)$, with possibly vanishing
	complex weights.
	The following conditions are equivalent:
	\begin{enumerate}
		\item[\textup{(i)}]
		The pair $\mathcal T$ supports a dense
		hypercyclic algebra.
		
		\item[\textup{(ii)}]
		The pair $\mathcal T$ supports a dense hypercyclic
		algebra freely generated, as a nonunital associative
		algebra, by countably infinitely many elements.
		
	\item[\textup{(iii)}]
	There exists a sequence of finite paths $(\sigma_k)$,
	with displacements $(m_k,n_k)$ satisfying
	$n_k\to\infty$, such that
	\[
	\min_{(i,j)\in F_N}
	|\Lambda_{\sigma_k}(i,j)|
	\longrightarrow\infty
	\qquad\text{for every }N\ge0.
	\]
	\end{enumerate}
\end{proposition}

\begin{proof}
	We first verify submultiplicativity.
	For $x,y\in\mathcal L_p$, put
	\[
	R_i=\left(\sum_{t\ge0}|x_{i,t}|^p\right)^{1/p},
	\qquad
	C_s=\left(\sum_{a=0}^{s}|y_{a,s-a}|^p\right)^{1/p}.
	\]
	Since $p\le2$, Cauchy--Schwarz and the inequality
	$\|\cdot\|_2\le\|\cdot\|_p$ give
	$|(x\ast y)_{i,j}|\le R_iC_{i+j}$. Hence
	\[
	\|x\ast y\|_p^p
	\le
	\sum_{i\ge0}R_i^p\sum_{s\ge i}C_s^p
	\le
	\|x\|_p^p\|y\|_p^p.
	\]
	Thus $(\mathcal L_p,\ast)$ is a separable Banach algebra.
	Clearly, \textup{(ii)} implies \textup{(i)}.
	
	\textup{(i)} $\Rightarrow$ \textup{(iii)}.
	Let $\mathcal A$ be a dense hypercyclic algebra.
	For each $k\ge1$, the commutator construction in the
	proof of Theorem~\ref{thm:triangular-shift-algebras},
	followed by scaling, gives
	\[
	0\ne z_k\in\mathcal A\cap J_k,
	\qquad
	\|z_k\|_p<\frac1{2k}.
	\]
	Choose a path $\sigma_k$ with displacement $(m_k,n_k)$
	such that
	\[
	\biggl\|
	T^{\sigma_k}z_k-
	\sum_{(i,j)\in F_k}e_{i,j}
	\biggr\|_p<\tfrac12.
	\]
	Its $(0,0)$ coordinate is nonzero.
	Since $z_k\in J_k$, this forces $n_k\ge k$.
	For every $(i,j)\in F_k$,
	\[
	\tfrac12
	<
	\left|
	\Lambda_{\sigma_k}(i,j)
	(z_k)_{i+m_k,j+n_k}
	\right|
	\le
	|\Lambda_{\sigma_k}(i,j)|\,\|z_k\|_p.
	\]
	Consequently,
	\[
	\min_{(i,j)\in F_k}
	|\Lambda_{\sigma_k}(i,j)|>k.
	\]
	For every fixed $N\ge0$ and $k\ge\max\{1,N\}$,
	the inclusion $F_N\subseteq F_k$ yields
	\[
	\min_{(i,j)\in F_N}
	|\Lambda_{\sigma_k}(i,j)|
	\ge
	\min_{(i,j)\in F_k}
	|\Lambda_{\sigma_k}(i,j)|
	>k.
	\]
	This proves \textup{(iii)}.
	
	\textup{(iii)} $\Rightarrow$ \textup{(ii)}.
	Fix $d\ge1$ and a finite nonempty set $A$ of nonempty
	words in $d$ variables. Choose a longest word
	$\beta\in A$ of length $\ell$.
	Given nonempty open sets $U_1,\ldots,U_d,V$, choose
	finitely supported $x_a\in U_a$ and $y\in V$,
	with all supports contained in $F_N$, and
	$\varepsilon>0$ such that
	$B(x_a,\varepsilon)\subseteq U_a$.
	
	By \textup{(iii)}, choose $k$ sufficiently large that
	\[
	n_k>2N+(N+1)^2(\ell-1)
	\]
	and $\min_{(i,j)\in F_N}|\Lambda_{\sigma_k}(i,j)|$
	is sufficiently large.
	Set $\sigma=\sigma_k$ and use the matrix-unit
	construction from
	Theorem~\ref{thm:triangular-shift-algebras}.
	Give every entry in the sequence indexed by $(i,j)$
	the same coefficient $\rho_{i,j}$, where
	\[
	\rho_{i,j}^{\ell}
	=
	\frac{y_{i,j}}{\Lambda_\sigma(i,j)}.
	\]
	The resulting finite perturbations satisfy
	\[
	\|u_a-x_a\|_p
	\le
	\|u_a-x_a\|_1
	\le
	\ell\sum_{(i,j)\in F_N}
	\left|
	\frac{y_{i,j}}{\Lambda_\sigma(i,j)}
	\right|^{1/\ell}
	<\varepsilon.
	\]
	Thus $u_a\in U_a$.
	The support and index calculation proving
	\eqref{eq:triangular-word-separation} is unchanged.
	Products along complete sequences equal
	$\rho_{i,j}^{\ell}$, so the same global identities hold:
	\[
	T^\sigma\beta(\mathbf u)=y,
	\qquad
	T^\sigma w(\mathbf u)=0
	\quad(w\in A\setminus\{\beta\}).
	\]
	Since $\beta$ depends only on $A$,
	Theorem~\ref{thm:algebra-transitivity-criterion}
	applies for every $d\ge1$ and proves \textup{(ii)}.
\end{proof}
\begin{remark}\label{rem:l2-triangular-algebras}
	The restriction $p\le2$ is essential for continuity
	of the triangular product. For $p>2$, put
	\[
	x_N=N^{-1/p}\sum_{t=1}^{N}e_{0,t},
	\qquad
	y_N=N^{-1/p}\sum_{t=1}^{N}e_{t,N+1-t}.
	\]
	Then $\|x_N\|_p=\|y_N\|_p=1$, whereas
	\[
	x_N\ast y_N=N^{1-2/p}e_{0,N+1},
	\qquad
	\|x_N\ast y_N\|_p\longrightarrow\infty.
	\]
	Thus the triangular product does not make
	$\mathcal L_p$ a Banach algebra with its usual norm
	when $p>2$.
\end{remark}
\providecommand{\smallodot}{\mathbin{\scriptstyle\odot}}

	The existence of a dense hypercyclic algebra can depend on
	the multiplication, even for a fixed operator pair on the
	same Fr\'echet space. The following example compares
	coordinatewise and triangular multiplication.
	
	\begin{example}\label{ex:product-dependent-algebras}
		On $X=\omega^2$ with the product topology, define
		\[
		(T_1x)_{i,j}=x_{i+1,j},\qquad
		(T_2x)_{i,j}=
		\begin{cases}
			(i+1)x_{i,1},&j=0,\\
			0,&j\ge1.
		\end{cases}
		\]
		The noncommuting pair $\mathcal T=(T_1,T_2)$ has the
		following properties on this same Fr\'echet space:
		\begin{enumerate}
			\item[\textup{(i)}]
			With coordinatewise multiplication, it supports
			a dense hypercyclic algebra.
			
			\item[\textup{(ii)}]
			With the triangular product, it supports a nonzero,
			singly generated hypercyclic algebra, but no dense
			hypercyclic algebra.
		\end{enumerate}
	\end{example}
	
	\begin{proof}
		Both operators are continuous, and
		\[
		T_1T_2e_{1,1}=2e_{0,0},\qquad
		T_2T_1e_{1,1}=e_{0,0},
		\]
		so they do not commute.
		Since all effective weights of $T_1$ equal $1$,
		Theorem~\ref{prop:omega2-pointwise-algebras}
		gives \textup{(i)}.
		
		For the triangular product, take $N=L=1$ in the path
		condition of Theorem~\ref{thm:triangular-shift-algebras}.
		Every path with $n\ge1$ has zero weight at $(0,1)$:
		its first second-coordinate step encounters
		$w^{(2)}_{a,2}=0$ for some $a\ge0$.
		The path condition therefore fails, so $\mathcal T$
		supports no dense hypercyclic algebra.
		
		It remains to construct a nonzero hypercyclic algebra.
		The shift $T_1$ is mixing, since restrictions on $F_N$
		and on $F_N+(k,0)$ can be prescribed independently
		for $k>N$. Moreover, it is an algebra endomorphism
		for the triangular product:
		\[
		T_1(u\ast v)=(T_1u)\ast(T_1v).
		\]
		For every nonzero $P\in\mathbb C[t]$ with $P(0)=0$,
		we claim that the continuous map
		\[
		\Phi_P:X\longrightarrow X,\qquad
		\Phi_P(z)=P_\ast(z),
		\]
		has dense range.
		
		Fix $M\ge0$ and let $B$ be an upper triangular
		$(M+1)\times(M+1)$ matrix with distinct diagonal
		entries $\mu_0,\ldots,\mu_M$.
		Choose $\alpha_s\in\mathbb C$ satisfying
		$P(\alpha_s)=\mu_s$, and interpolate a polynomial
		$Q$ such that $Q(\mu_s)=\alpha_s$.
		The polynomial $P(Q(t))-t$ vanishes at every $\mu_s$,
		so it is divisible by the characteristic polynomial
		of $B$. The Cayley--Hamilton theorem gives
		\[
		P(Q(B))=B.
		\]
		Since $\pi_M$ is a surjective algebra homomorphism,
		it follows that $B\in\pi_M(\Phi_P(X))$.
		Matrices with distinct diagonal entries are dense
		in the finite upper triangular matrix algebra.
		Thus $\pi_M(\Phi_P(X))$ is dense for every $M$.
		As these blocks capture every finite coordinate set,
		$\Phi_P(X)$ is dense in $X$.
		
		Choose $x\in HC(T_1)$. For every such polynomial $P$,
		the endomorphism property yields
		\[
		T_1^kP_\ast(x)=\Phi_P(T_1^kx)\qquad(k\ge0).
		\]
		Continuity and density imply
		\[
		\overline{\{T_1^kP_\ast(x):k\ge0\}}
		\supseteq\overline{\Phi_P(X)}=X.
		\]
		Hence $P_\ast(x)$ is nonzero and hypercyclic for $T_1$,
		and therefore for $\mathcal T$.
		Consequently,
		\[
		\operatorname{alg}_\ast(x)
		=
		\{P_\ast(x):P\in\mathbb C[t],\ P(0)=0\}
		\]
		is a nonzero, singly generated hypercyclic algebra.
		This proves \textup{(ii)}.
	\end{proof}
	
	\begin{remark}\label{rem:density-essential}
		Example~\ref{ex:product-dependent-algebras} shows that
		a nonzero hypercyclic algebra need not imply the existence
		of a dense one for the triangular product.
		The same generator satisfies
		$x^{\ast r}\in HC(\mathcal T)$ for every $r\ge1$.
		Therefore, the restriction $x\in J_1$ in
		Theorem~\ref{thm:triangular-shift-algebras}\textup{(i)}
		cannot be omitted.
	\end{remark}
	
	\backmatter
	
	\section*{Declarations}
	
	\bmhead{Author Contributions}
	Xiang Chen made the primary contribution to the research and
	prepared the initial draft.  All authors read
	and approved the final manuscript.
	
	\bmhead{Conflict of Interest}
	The authors declare no competing interests with the journal's editors, reviewers, or readers.
	
	\bmhead{Funding}
	This work was supported by  the National Natural Science Foundation of
	China (Grant Nos. 12571088).
	
	\bmhead{Data Availability Statement}
	No datasets were generated or analyzed in this theoretical study.

\end{document}